\documentclass[12pt]{amsart}
\usepackage[left=1in, right=1in, top=1in, bottom=1in]{geometry}
\usepackage[margin=20pt]{caption} 
\usepackage{framed}

\usepackage{amsmath}
\usepackage{amssymb}
\usepackage{amsthm}
\usepackage{babel}
\usepackage{mathrsfs}
\usepackage{hyperref}
\usepackage{graphicx}
\usepackage{wrapfig}
\usepackage{comment}
\usepackage{xcolor}
\usepackage{comment}
\usepackage{import}
\usepackage{calc}
\usepackage{enumitem}

\usepackage{caption}

\usepackage{multirow} 

\newtheorem{theorem}{Theorem}

\newtheorem*{obs*}{Observation}

\theoremstyle{definition}
\newtheorem{definition}{Definition}

\newcommand{\Z}{\mathbb{Z}}

\newcommand{\Q}{\mathbb{Q}}

\DeclareMathOperator{\h}{H} 
\DeclareMathOperator{\lk}{lk}

\newcommand{\clasper}{\textsf{Clasper}} 

\begin{document}
\title{An Algorithm to calculate Generalized Seifert Matrices}

\author{Stefan Friedl}
\address{Stefan Friedl, Universit\"at Regensburg, 93047 Regensburg, Germany}
\email{sfriedl@gmail.com}

\author{Chinmaya Kausik}
\address{Chinmaya Kausik, Department of Mathematics, University of Michigan, Ann Arbor, USA; 48104}
\email{chinmaya.kausik.1@gmail.com}

\author{Jos\'e Pedro Quintanilha}
\address{Fakultät für Mathematik, Universität Bielefeld, Postfach 100131, D-33501 Bielefeld, Germany}
\email{zepedro.quintanilha@gmail.com}

\begin{abstract}
We develop an algorithm for computing generalized Seifert matrices for colored links given as closures of colored braids. The algorithm has been implemented by the second author as a computer program called \clasper. \clasper\ also outputs the Conway potential function, the multivariable Alexander polynomial and the Cimasoni-Florens signatures of a link, and displays a visualization of the C-complex used for producing the generalized Seifert matrices.
\end{abstract}

\maketitle

\section{Introduction}

\subsection{Background}
A Seifert surface for an oriented knot $K\subset S^3$ is a connected compact oriented surface $F\subset S^3$ with $\partial F=K$. Recall  \cite[Chapter~VIII]{Ro90} that the corresponding Seifert form is defined as
\[ \begin{array}{rcl} \h_1(S)\times \h_1(S)&\to & \Z\\
([\gamma],[\delta])&\mapsto & \operatorname{lk}(\gamma^+,\delta),\end{array}\]
where $\gamma^+$ denotes a positive push-off of $\gamma$ and $\operatorname{lk}$ denotes the linking number of oriented disjoint curves in $S^3$. Any choice of a basis for $\h_1(S)$
now defines a corresponding matrix, called Seifert matrix.  Seifert matrices,
have played an important role in knot theory
ever since they were introduced by Herbert Seifert \cite{Se34}.
For example a Seifert matrix~$A$ can be used to calculate the Alexander polynomial $\Delta_K(t)=\operatorname{det}(At-A^T)\in \mathbb{Z}[t^{\pm 1}]$ and
it can be used to define the
 Levine-Tristram signature function $\sigma(K)\colon S^1\to \mathbb{Z}$
by setting $\sigma_z(K):=\operatorname{sign}(A(1-z)+A^T(1-\overline{z}))$. 
There are various algorithms for computing Seifert matrices for a knot 
\cite{O'B02,Col16}.  In particular, Julia Collins gave an algorithm to determine a Seifert matrix from a given braid description \cite{Col16}. An implementation of this algorithm is available online \cite{CKL16}.

\begin{figure}[h]
\centering
\def\svgwidth{0.45 \linewidth}
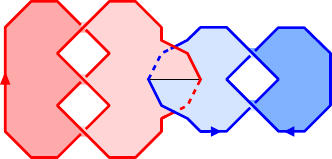
 \caption{A C-complex for a $3$-component link in $\mu = 2$ colors. The Seifert surface for the trefoil knot (in red) and the Seifert surface for the Hopf link (in blue) intersect along a clasp in the center of the picture.}
\label{fig.clasp_complex}
\end{figure}

Less known are the generalizations of Seifert surfaces and Seifert matrices to links.
Let $L=L_1\sqcup \dots\sqcup L_\mu\subset S^3$ be a $\mu$-colored oriented link, i.e.\ $L$ is a disjoint union of finitely many oriented knots that get grouped into $\mu$ sets.
Daryl Cooper \cite{Coo82} and David Cimasoni \cite{Ci04} introduced the notion of a C-complex for $L$. A C-complex consists, roughly speaking, of $\mu$ embedded compact oriented surfaces $S_1\cup \dots\cup S_\mu$ with $\partial S_i=L_i$ and a few restrictions on how the $S_i$ are allowed to intersect. We postpone the definition to Section~\ref{section:c-complex}, but we hope that Figure~\ref{fig.clasp_complex} gives at least an idea of the concept.

Given a C-complex $S$ and given a basis for $\h_1(S)$
one obtains for any choice of $\epsilon\in \{\pm 1\}^\mu$ a  generalized Seifert matrix $A^\epsilon$. 
These  $2^\mu$ generalized Seifert matrices 
can be used to define and calculate the Conway potential function, which in turn determines the multivariable Alexander polynomial \cite[p.~128]{Ci04} \cite{DMO21}. We recall the formulae in 
Theorem~\ref{thm:conway-potential-function}.
Furthermore, David Cimasoni and Vincent Florens used generalized Seifert matrices to define a generalization of the Levine-Tristram signature function \cite{CF08}, namely the Cimasoni-Florens signature function $\sigma_L\colon (S^1)^\mu\to \Z$.
The generalized Seifert matrices can also be used to determine the Blanchfield pairing of a colored link \cite{CFT18,Con18}.

\subsection{Our results}

Our goal was to come up with an algorithm that computes generalized Seifert matrices of colored links, and implement it as a computer program. To formulate an algorithm one first has to settle on the input.
The usual proof of Alexander's Theorem (see for example the account of Burde-Zieschang-Heusener \cite[Proposition~2.12]{BZH14}) shows that every oriented $\mu$-colored link is the closure of a $\mu$-colored braid (i.e.\  a braid together with an integer $0\le k < \mu$ associated to each component of its closure). The description of braids as sequences of elementary crossings makes them a convenient input type for a computer program.

\begin{framed}
We give an  algorithm that takes as input a colored braid, produces a C-complex~$S$ for its closure, and computes the associated generalized Seifert matrices $A^\epsilon$ (with respect to some basis of~$\h_1(S)$).
\end{framed}

In this paper, our algorithm is explained in natural language. Despite its geometrical flavor, it is formal enough to be implemented as a computer program:

\begin{framed}
We provide a computer program, called \clasper, implementing our algorithm. \clasper\ displays a visualization of the constructed C-complex and outputs a family of generalized Seifert matrices. It  also computes the Conway potential function, the multivariable Alexander polynomial, and the Cimasoni-Florens signatures of the colored link.
\end{framed}

\begin{figure}[h]
    \centering
    \includegraphics[width=0.95 \linewidth]{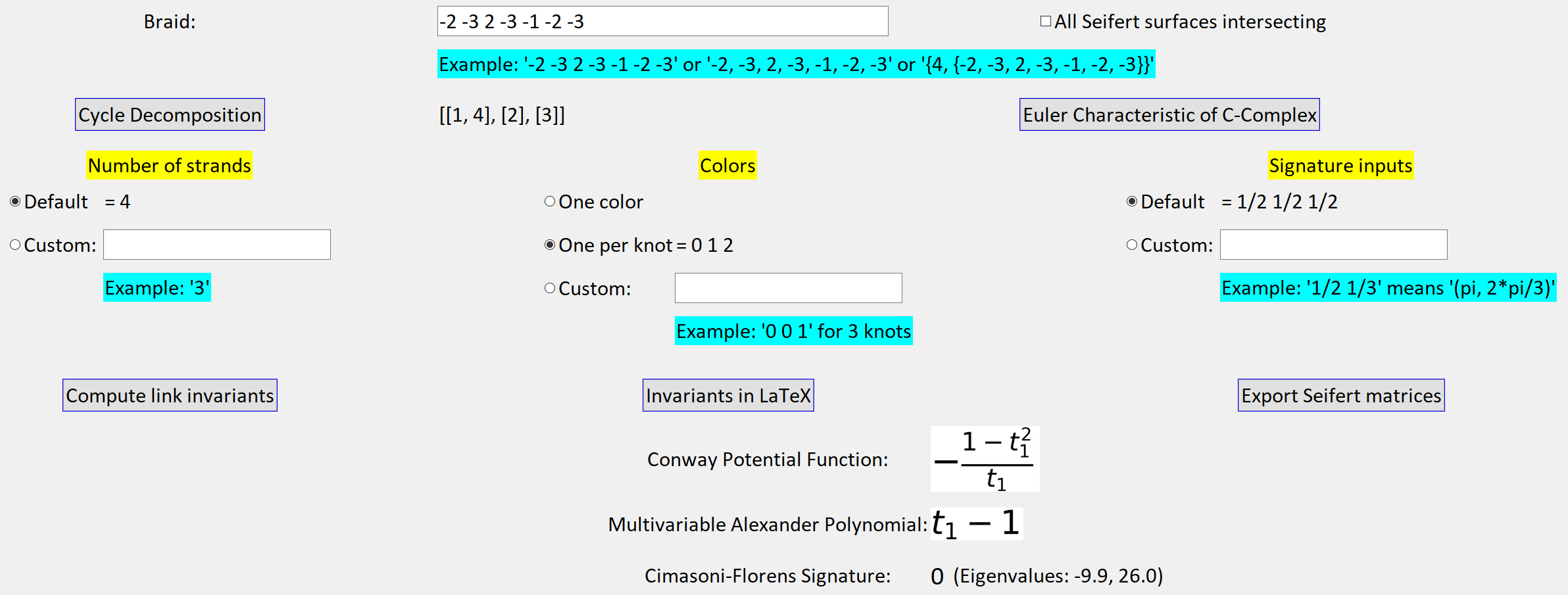}\medskip

    \includegraphics[width=0.85 \linewidth]{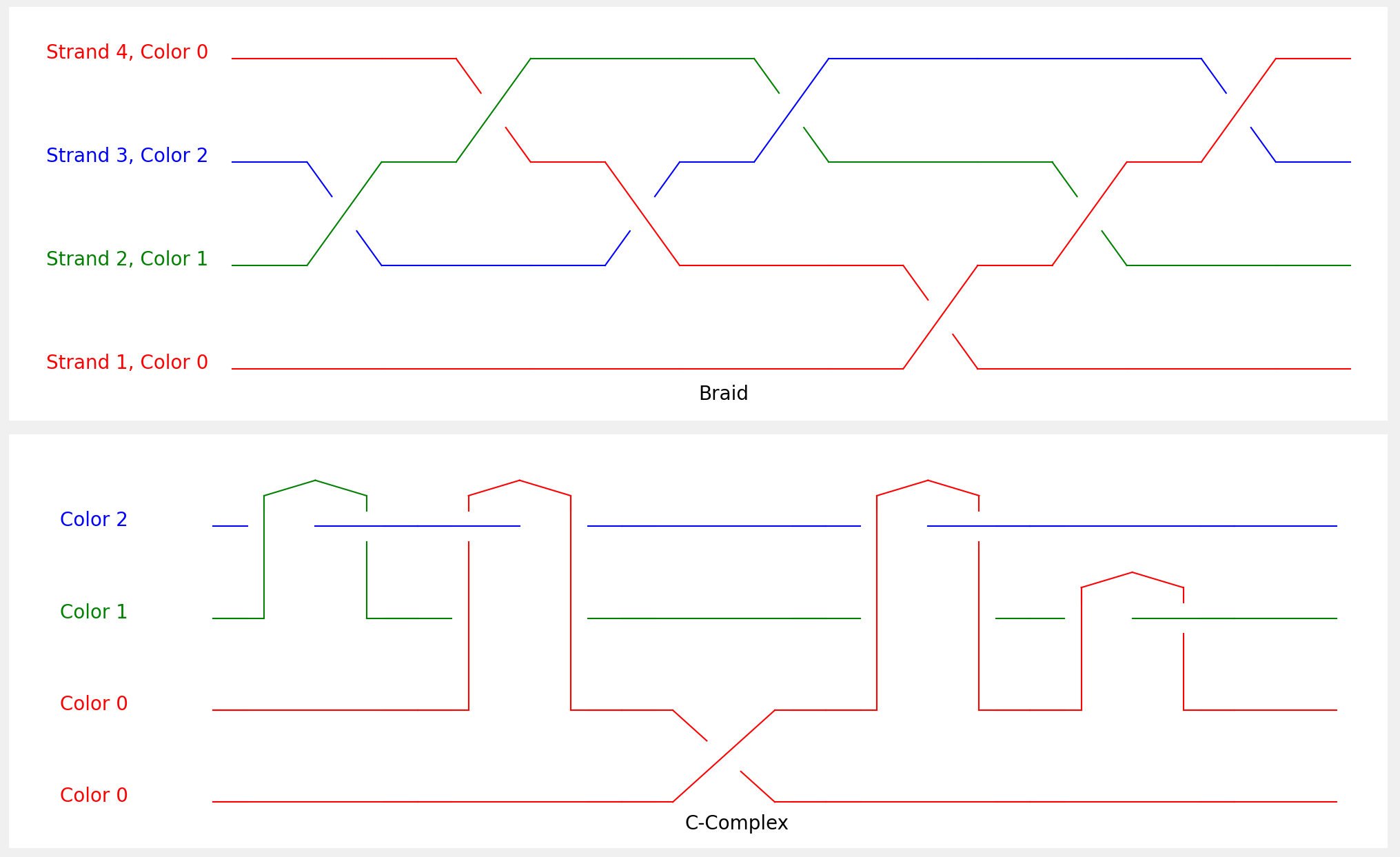}
    \caption{The \clasper\ interface. Top: \clasper\ takes as input a braid given as its sequence of crossings, number of strands, and information about the coloring of its closure, and outputs a family of generalized Seifert matrices thereof. \clasper\ also computes the Conway potential function, multivariable Alexander polynomial, and Cimasoni-Florens signature (at a user-specified point). Bottom: \clasper\ displays a diagram of the input colored braid and a schematic of the C-complex from which the generalized Seifert matrices were produced. See Subsection~\ref{sec.fillin} for how to interpret such a schematic, and Figure~\ref{fig.multicolorspine} (right) for the C-complex represented in this screenshot.}
    \label{fig.screenshots}
\end{figure}

\clasper\ was programmed by the second author. A Windows installer, as well as the Python source code, can be downloaded at \url{https://github.com/Chinmaya-Kausik/py_knots}. Figure~\ref{fig.screenshots} shows the user interface of \clasper. 

\subsection{Organization of the article}
In Section~\ref{section:c-complex-matrix} we give the definition of C-complexes and we show how they can be used to define generalized Seifert matrices. Furthermore we recall how generalized Seifert matrices can be used to compute the multivariable Alexander polynomial and how they can be used to define the Cimasoni-Florens signatures.
In Section~\ref{section:algorithm} we explain our algorithm for computing generalized Seifert matrices for colored links given by a braid description. 
Section~\ref{section:implementation} contains additional remarks on the technical details of the implementation by the second author.

\subsection*{Acknowledgments} 
SF and JPQ were supported by the SFB 1085 ``higher invariants'' at the University of Regensburg, funded by the DFG. We also wish to thank Lars Munser for helpful conversations.

\section{C-complexes and generalized Seifert matrices}\label{section:c-complex-matrix}

\subsection{C-complexes}\label{section:c-complex}

\begin{definition}
A \textbf{C-complex} (where ``C'' is short for ``clasp'') for a colored  oriented link $L =L_1\cup \dots \cup L_\mu\subset S^3$ is a collection of  surfaces $S_1,\dots,S_\mu$ such that:
\begin{itemize}
    \item  each $S_i$ is a compact oriented surface in $S^3$ with $\partial S_i=L_i$, where we demand that the equality is an equality of oriented 1-manifolds,
  \item every two distinct $S_i, S_j$ intersect transversely along a (possibly empty) finite disjoint union of intervals, each with one endpoint in~$L_i$ and the other in~$L_j$, but otherwise disjoint from~$L$. Such an intersection component is called a \textbf{clasp} -- see Figure~\ref{fig.clasp_complex} for an illustration, 
\item the union $S_1\cup\dots\cup S_\mu$ is connected, and
    \item there are no triple points: for distinct $i,j,k$, we have $S_i \cap S_j \cap S_k = \emptyset$.
\end{itemize}
\end{definition}

C-complexes were introduced for 2-component links by Cooper \cite{Coo82}
and in the general case by Cimasoni \cite{Ci04}.
Cimasoni \cite[Lemma~1]{Ci04} showed that every colored oriented link admits a C-complex.
In this paper we will give a different proof, in fact we will outline an algorithm which takes as input a braid description of a colored oriented link and which produces an explicit C-complex.

We remark that any two C-complexes for a colored link are related by a finite sequence of moves of certain types, a complete list of which was given by by Davis-Martin-Otto \cite[Theorem~1.3]{DMO21} building on work of Cimasoni \cite{Ci04}. We will however not make use of this fact.

\subsection{Push-offs of curves and generalized Seifert matrices}\label{sec.seifertmatrices}
Let $L$ be a $\mu$-colored oriented link  and let $S_1,\dots,S_\mu$ be a C-complex for $L$.
Following Cooper \cite{Coo82} and Cimasoni \cite{Ci04} we will now associate to this data generalized Seifert pairings and generalized Seifert matrices.  The approach to defining generalized Seifert pairings and matrices is quite similar to the more familiar definition for knots recalled in the introduction. 

\begin{figure}[h]
    \centering
    \def\svgwidth{0.8\linewidth}
    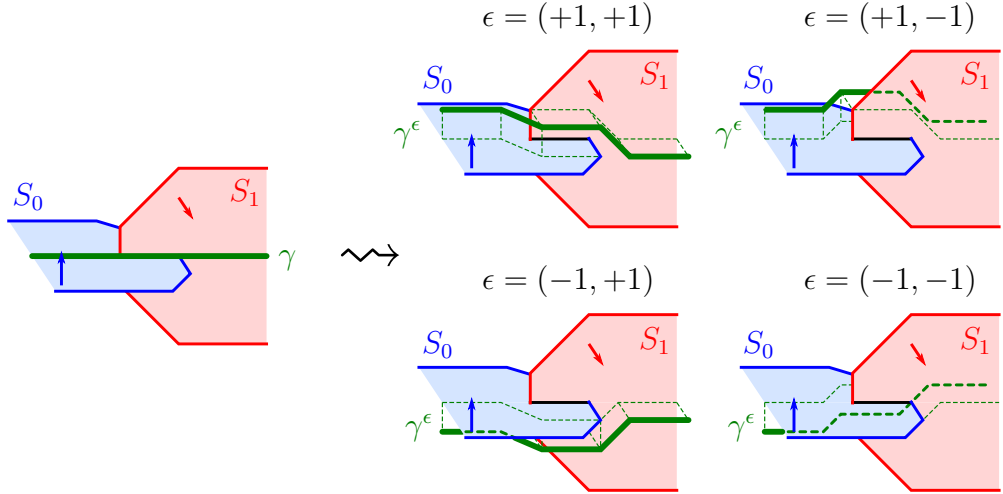
    \caption{Push-offs of a path~$\gamma$ in a C-complex determined by a $\mu$-tuple~$\epsilon$. Arrows indicate the framings of the normal bundles of the Seifert surfaces induced by their orientations. In this example, the C-complex is comprised of $\mu = 2$ Seifert surfaces. We show all its push-offs~$\gamma^\epsilon$, determined by the four possible pairs $\epsilon \in \{\pm 1\}^\mu$.}
    \label{fig.clasp_pushoff}
\end{figure}

First note that the orientation of the link $L$ induces  an orientation on  each Seifert surface~$S_i$ in our C-complex~$S$, which in turn induces an orientation of the normal bundle of~$S_i$. Now, each of the $2^\mu$ tuples of signs $\epsilon = (\epsilon_1, \ldots, \epsilon_\mu)$, with $\epsilon_i = \pm 1$, prescribes a way of pushing a point $p\in S$ off of~$S$: at the Seifert surface $S_i$ and away from the clasps, we let the sign~$\epsilon_i$ determine whether to push $p$ in the positive or negative direction of the normal bundle of~$S_i$, and at a clasp between $S_i$~and~$S_j$, we move~$p$ in the ``diagonal'' direction specified by $\epsilon_i, \epsilon_j$. The push-off of a path~$\gamma$ in~$S$ specified by a tuple $\epsilon  \in \{\pm 1\}^\mu$ will be denoted by $\gamma^\epsilon$; see Figure~\ref{fig.clasp_pushoff}. 

Given any $\epsilon\in \{\pm 1\}^\mu$ Cimasoni defines the generalized Seifert pairing 
\[ \begin{array}{rcl}\alpha^\epsilon\colon \h_1(S) \times \h_1(S)& \to &\Z\\
([\gamma],[\delta])&\mapsto & \operatorname{lk}(\gamma^\epsilon, \delta).\end{array}\]
Finally we pick  a basis of $\h_1(S)$. 
The collection of matrices  $A^\epsilon$ of $\alpha^\epsilon$ with respect to this basis is called a collection of \textbf{generalized Seifert matrices}  of $L$.

\subsection{The Conway potential function, the multivariable Alexander polynomial and  generalized signatures}\label{sec.invariants}
In this section we turn to the discussion of several applications of generalized Seifert matrices.

In 1970 John Conway \cite{Con70} associated to any $\mu$-colored link~$L$ a rational function on $\mu$ variables, now called the Conway potential function $\nabla_L(t_1,\dots,t_\mu)\in \Q(t_1,\dots,t_\mu)$.
Cimasoni showed that the Conway potential function can be computed using generalized Seifert matrices. To state Cimasoni's Theorem we need the following definition. 

\begin{figure}[h]
\centering
\def\svgwidth{0.45 \linewidth}
\begingroup%
  \makeatletter%
  \providecommand\color[2][]{%
    \errmessage{(Inkscape) Color is used for the text in Inkscape, but the package 'color.sty' is not loaded}%
    \renewcommand\color[2][]{}%
  }%
  \providecommand\transparent[1]{%
    \errmessage{(Inkscape) Transparency is used (non-zero) for the text in Inkscape, but the package 'transparent.sty' is not loaded}%
    \renewcommand\transparent[1]{}%
  }%
  \providecommand\rotatebox[2]{#2}%
  \newcommand*\fsize{\dimexpr\f@size pt\relax}%
  \newcommand*\lineheight[1]{\fontsize{\fsize}{#1\fsize}\selectfont}%
  \ifx\svgwidth\undefined%
    \setlength{\unitlength}{71.99999892bp}%
    \ifx\svgscale\undefined%
      \relax%
    \else%
      \setlength{\unitlength}{\unitlength * \real{\svgscale}}%
    \fi%
  \else%
    \setlength{\unitlength}{\svgwidth}%
  \fi%
  \global\let\svgwidth\undefined%
  \global\let\svgscale\undefined%
  \makeatother%
  \begin{picture}(1,0.34409576)%
    \lineheight{1}%
    \setlength\tabcolsep{0pt}%
    \put(0,0){\includegraphics[width=\unitlength,page=1]{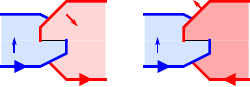}}%
    \put(0.20122718,0.20172868){\makebox(0,0)[lt]{\lineheight{1.25}\smash{\begin{tabular}[t]{l}$C$\end{tabular}}}}%
    \put(0.09783279,0.16810272){\color[rgb]{0,0.50196078,0}\makebox(0,0)[lt]{\lineheight{1.25}\smash{\begin{tabular}[t]{l}$P$\end{tabular}}}}%
    \put(0,0){\includegraphics[width=\unitlength,page=2]{clasp_sign.pdf}}%
    \put(0.77414384,0.20172868){\makebox(0,0)[lt]{\lineheight{1.25}\smash{\begin{tabular}[t]{l}$C$\end{tabular}}}}%
    \put(0.67074962,0.16810272){\color[rgb]{0,0.50196078,0}\makebox(0,0)[lt]{\lineheight{1.25}\smash{\begin{tabular}[t]{l}$P$\end{tabular}}}}%
    \put(0,0){\includegraphics[width=\unitlength,page=3]{clasp_sign.pdf}}%
  \end{picture}%
\endgroup%

 \caption{The sign of a clasp~$C$. Link orientations and induced normal framings on the Seifert surfaces are indicated by arrows. Left: a negative, or left-handed clasp. Right: a positive, or right-handed clasp.}
\label{fig.clasp_sign}
\end{figure}

\begin{definition}
Let $S=(S_1,\dots,S_\mu)$ be a C-complex for a $\mu$-colored link. We define the \textbf{sign of a clasp}~$C$ by choosing one of its endpoints~$P$, and then taking the sign of the intersection between the Seifert surface and the link component at~$P$. This definition is independent of the choice of endpoint~$P$; see Figure~\ref{fig.clasp_sign}.

The \textbf{sign of a C-complex}~$S$ is defined to be the product of the signs of all clasps, and denoted by~$\operatorname{sgn}(S)$.
\end{definition}

The following theorem  of Cimasoni \cite[p.~128]{Ci04} gives a way to calculate the Conway potential function.

\begin{theorem}\label{thm:conway-potential-function}
Let $L$ be a $\mu$-colored oriented link and let $S$ be a C-complex for $L$.
Choose a basis for $\h_1(S)$ and use it to define the $2^\mu$ generalized  Seifert matrices
$A^\epsilon$ with $\epsilon=(\epsilon_1,\dots,\epsilon_\mu)\in\{\pm 1\}^\mu$.
Then
\[\nabla_L(t_1,\dots,t_\mu) = \operatorname{sgn}(S)\cdot \prod\limits_{i=1}^\mu
(t_i-t_i^{-1})^{-1+\chi\big(\bigcup_{j\ne i} S_j \big)}\cdot \operatorname{det}\Big(-\sum_{\epsilon \in \{\pm 1\} ^\mu} A^\epsilon\cdot \epsilon_1 \dots \epsilon_\mu\cdot  t_1^{\epsilon_1}\dots t_\mu^{\epsilon_\mu}\Big).\]
\end{theorem}

In fact the right-hand side of Theorem~\ref{thm:conway-potential-function} can be used as a definition of the Conway potential function \cite[Lemma~4]{Ci04} \cite[Lemma~2.1]{DMO21}. Also note that for $\mu=1$ the above definition  differs from the definition of the one-variable Conway polynomial~$\nabla_L(z)$
given by Lickorish \cite{Li97} and LinkInfo \cite{LM22} by the substitution $z=t-t^{-1}$.

Next recall that given a $\mu$-colored oriented link $L$  we can use presentation matrices for the Alexander module $\h_1(S^3\setminus L;\Z[t_1^{\pm 1},\dots,t_\mu^{\pm 1}])$
to define the multivariable Alexander polynomial $\Delta_L(t_1,\dots,t_\mu)$, which is well-defined only up to multiplication by a term of the form $\pm t_1^{k_1}\dots t_\mu^{k_\mu}$, with $k_i \in \Z$.
The Conway potential function can be viewed as a refinement of the multivariable Alexander polynomial $\Delta_L(t_1,\dots,t_\mu)$. More precisely Conway shows \cite[p.~338]{Con70} that, up to the above indeterminacy, the following equality holds
\[ \nabla_L(t_1,\dots,t_\mu)\,\,=\,\, \left\{ \begin{array}{ll}
\frac{1}{t_1-t_1^{-1}}\cdot \Delta_L(t_1^2), &\mbox{ if }\mu=1,\\
\Delta_L(t_1^2,\dots,t_\mu^2), &\mbox{ if }\mu \geq 2.\end{array}\right.\]
It follows from this equality that the multivariable Alexander polynomial $\Delta_L(t_1,\dots,t_\mu)$ is determined by the Conway potential function $\nabla_L(t_1,\dots,t_\mu)$. In particular, in light of Theorem~\ref{thm:conway-potential-function},  it can be determined by the generalized Seifert matrices. 

\begin{definition}\label{dfn.CFsignatures}
Let $L$ be a $\mu$-colored oriented link and let $S$ be a C-complex for $L$.
We pick a basis for $\h_1(S)$ and we use it to define the generalized Seifert matrices
$A^\epsilon$. Following Cimasoni-Florens \cite{CF08} we define
\begin{center}
$\displaystyle H(\omega) := \prod_{i=1}^\mu (1 - \overline{\omega}_i) \cdot A(\omega_1,\dots,\omega_\mu)$,
\end{center}
where
\[A(t_1, \ldots , t_\mu) := \sum_{\epsilon \in \{0,1\}^\mu}  \epsilon_1 \cdots \epsilon_\mu \cdot t_1^\frac{1-\epsilon_1}{2} \cdots  t_\mu^\frac{1-\epsilon_\mu}{2} \cdot A^\epsilon\]
and
$\omega = (\omega_1,\dots,\omega_\mu) \in (S^1 \setminus \{1\})^\mu$. We define the \textbf{generalized signature} of $L$ at $\omega$ as 
\[ \sigma_L(\omega)\,\,=\,\, \mbox{signature of the hermitian matrix $H(\omega)$}\]
and we define the \textbf{nullity} of $L$ at $\omega$ as 
\[ \eta_L(\omega)\,\,=\,\,b_0(S)-1+ \mbox{nullity of the matrix $H(\omega)$}.\]
\end{definition}

We have the following theorem due to Cimasoni--Florens \cite[Theorem~2.1]{CF08}
and Davis--Martin--Otto \cite[Theorem~3.2]{DMO21}.

\begin{theorem}
Let  $\omega = (\omega_1,\dots,\omega_\mu) \in (S^1 \setminus \{1\})^\mu$.
The signature and nullity of a $\mu$-colored oriented link at $\omega$ are  invariant under isotopy of colored links. 
\end{theorem}

As is explained in \cite{CF08}, the signature and nullity invariants of links contain a lot of deep information on links:
\begin{enumerate}
    \item If $L$ is a $\mu$-colored oriented link that is isotopic to its mirror image,
then the signature function is identically zero  \cite[Corollary~2.11]{CF08}.
\item Two $m$-component oriented links $L=L_1\sqcup\dots\sqcup L_m$ and $J=J_1\sqcup\dots\sqcup J_m$ are called \textbf{smoothly concordant} if there exist disjoint properly smoothly embedded oriented annuli $A_1,\dots,A_m\subset S^3\times [0,1]$ such that 
$\partial A_i=L_i\times \{0\}\cup (-J_i)\times \{1\}$. If we treat $L$ and $J$ as $m$-colored links in the obvious way, then \cite[Theorem~7.1]{CF08}
for every prime power $q=p^k$ and every $\omega_1,\dots,\omega_m\in S^1\setminus \{1\}$ with $\omega_1^q=\dots=\omega_m^q=1$ we have 
\[ \sigma_L(\omega_1,\dots,\omega_m)\,=\, \sigma_J(\omega_1,\dots,\omega_m).\]
\end{enumerate}
We also point the reader towards two other applications:
\begin{enumerate}
    \item If $L$ is a $\mu$-colored link and $S_1,\dots,S_\mu$ is a C-complex, such that for any $i\ne j$ we have $S_i\cap S_j\ne \emptyset$, then the generalized Seifert matrices 
can be used to give an explicit presentation matrix for the multivariable
 Alexander module $\h_1(S^3\setminus L;\Lambda_\mu)$, where $\Lambda_\mu:=\Z[t_1^{\pm 1},\dots,t_\mu^{\pm 1},(1-t_1)^{-1},\dots,(1-t_\mu)^{-1}]$ \cite[Theorem~3.2]{CF08}. 
 \item Generalized signatures can be used to calculate Casson-Gordon invariants \cite{CG75,CG78} (which are a special case of Atiyah-Patodi-Singer invariants \cite{APS75}) of surgeries $M$ on links, corresponding to characters $\chi\colon \operatorname{H}_1(M)\to S^1$ \cite[Theorem~6.4]{CF08}.
\end{enumerate}

\section{Explanation of the algorithm}\label{section:algorithm}

We now describe our algorithm, which takes as input a colored braid, and produces:
\begin{enumerate}
    \item a C-complex for it, encoded as a graph with decorations (which we will call a ``decorated spine''), and
    \item a family of generalized Seifert matrices for that C-complex (with respect to some homology basis).
\end{enumerate}

With these matrices in hand, the computation of the Conway potential function and the Cimasoni-Florens signatures are conceptually straightforward using Theorem~\ref{thm:conway-potential-function} and Definition~\ref{dfn.CFsignatures}.

In Subsection~\ref{sec.warmup}, we warm up by explaining how to produce a Seifert surface for the closure of a braid on only one color, and how to and encode this Seifert surface as a decorated spine. This should help the reader familiarize themselves with our conventions. Subsection~\ref{sec.fullalg} lays out the full algorithm for braids on multiple colors, exemplifying the construction of the C-complex on a running example, and its encoding as a decorated spine. Finally, Subsection~\ref{sec.readmatrix} explains how to produce a homology basis for the C-complex and construct its associated generalized Seifert matrices. This is essentially an analysis of how to read off the relevant linking numbers from the decorated spines.

\subsection{Constructing Seifert surfaces}\label{sec.warmup}

We warm up to the construction end encoding of the C-complex with the case where there is a single color -- in other words, we construct a Seifert surface for a link given as the closure of a braid.

\begin{figure}[h]
    \centering
    \def\svgwidth{0.9 \linewidth}
    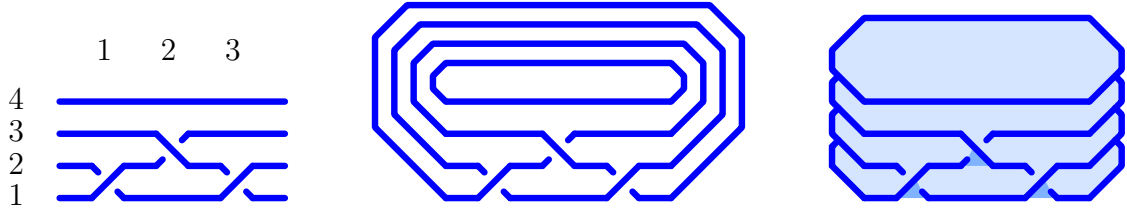
    \caption{Constructing a surface from a braid. Left: an example of a single-color braid. The numbers on the left indicate our convention for numbering the positions in a braid, and the numbers above the braid indicate the ordering of the crossings. Here, the number of strands is $n = 4$, and the sequence representing the braid is $s=(-1, 2, -1)$. Center: a diagram for the closure of the braid. Right: the surface obtained by applying Seifert's algorithm, visualized as a stack of Seifert disks connected by half-twisted bands. In this example, the surface is not connected.}
    \label{fig.1colorbraid}
\end{figure}

More concretely, the input data is a number~$n$ of strands, together with a sequence $s = (s_1, \ldots, s_m)$ of integers in~$\{-(n-1), \ldots, n -1 \}\backslash \{0\}$ specifying the crossings. We will adopt the convention of orienting the strands from left to right, and numbering the $n$~positions from bottom to top; see Figure \ref{fig.1colorbraid} (left) for an example. Each integer~$\pm s_i$ then represents a crossing between the strands in positions $s_i$~and~$s_i+1$, with a plus sign indicating that the over-crossing strand goes down one position (right-handed crossing), and a minus sign meaning that it goes up (left-handed crossing).

We close the braid by drawing $n$ arcs above it as illustrated in Figure~\ref{fig.1colorbraid} (center), and apply Seifert's algorithm to the resulting diagram. Explicitly, this means that each of the $n$ positions gives rise to a Seifert disk, and a crossing between the strands in positions $k$ and $k+1$ translates into a half-twisted band connecting the corresponding disks, the sign of the crossing determining the handedness of the twist. It is often convenient to visualize the Seifert disks as a ``stack of pancakes'', as in Figure~\ref{fig.1colorbraid} (right). Since all half-twisted bands connect Seifert disks in a manner that respects their top and bottom sides, the resulting surface is orientable.

\begin{figure}[h]
    \centering
    \def\svgwidth{0.4\linewidth}
    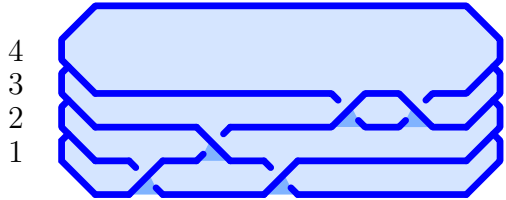
    \caption{Adding half-twisted bands to ensure connectedness. We obtain a Seifert surface from the surface in Figure~\ref{fig.1colorbraid} (right) by adding a pair of half-twisted bands with opposite handedness connecting the disks in positions $3$~and~$4$. This does not change the underlying link.}
    \label{fig.connected}
\end{figure}

The surface we produced might not yet be a Seifert surface, because it could be disconnected. We remedy this  by adding, for each pair of adjacent disks that are in different connected components, a pair of half-twisted bands with opposite handedness -- see Figure~\ref{fig.connected}. This is a harmless modification, as the diagram for the braid closure changes by a Reidemeister move of type II.

\begin{figure}[h]
    \centering
    \def\svgwidth{0.5\linewidth}
\begingroup%
  \makeatletter%
  \providecommand\color[2][]{%
    \errmessage{(Inkscape) Color is used for the text in Inkscape, but the package 'color.sty' is not loaded}%
    \renewcommand\color[2][]{}%
  }%
  \providecommand\transparent[1]{%
    \errmessage{(Inkscape) Transparency is used (non-zero) for the text in Inkscape, but the package 'transparent.sty' is not loaded}%
    \renewcommand\transparent[1]{}%
  }%
  \providecommand\rotatebox[2]{#2}%
  \newcommand*\fsize{\dimexpr\f@size pt\relax}%
  \newcommand*\lineheight[1]{\fontsize{\fsize}{#1\fsize}\selectfont}%
  \ifx\svgwidth\undefined%
    \setlength{\unitlength}{83.67548665bp}%
    \ifx\svgscale\undefined%
      \relax%
    \else%
      \setlength{\unitlength}{\unitlength * \real{\svgscale}}%
    \fi%
  \else%
    \setlength{\unitlength}{\svgwidth}%
  \fi%
  \global\let\svgwidth\undefined%
  \global\let\svgscale\undefined%
  \makeatother%
  \begin{picture}(1,0.42127057)%
    \lineheight{1}%
    \setlength\tabcolsep{0pt}%
    \put(0,0){\includegraphics[width=\unitlength,page=1]{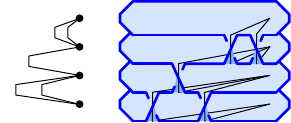}}%
    \put(0.00195918,0.09793724){\color[rgb]{0,0,0}\makebox(0,0)[lt]{\lineheight{1.25}\smash{\begin{tabular}[t]{l}$-$\end{tabular}}}}%
    \put(0.0472261,0.1988978){\color[rgb]{0,0,0}\makebox(0,0)[lt]{\lineheight{1.25}\smash{\begin{tabular}[t]{l}$+$\end{tabular}}}}%
    \put(0.15681555,0.09777674){\color[rgb]{0,0,0}\makebox(0,0)[lt]{\lineheight{1.25}\smash{\begin{tabular}[t]{l}$-$\end{tabular}}}}%
    \put(0.13572018,0.29564343){\color[rgb]{0,0,0}\makebox(0,0)[lt]{\lineheight{1.25}\smash{\begin{tabular}[t]{l}$-$\end{tabular}}}}%
    \put(0.25288298,0.29717376){\color[rgb]{0,0,0}\makebox(0,0)[lt]{\lineheight{1.25}\smash{\begin{tabular}[t]{l}$+$\end{tabular}}}}%
    \put(0,0){\includegraphics[width=\unitlength,page=2]{1colorspine.pdf}}%
  \end{picture}%
\endgroup%

        \caption{A decorated spine. Left: the spine for the Seifert surface in Figure~\ref{fig.connected}. Edges are labeled with signs indicating the handedness of the corresponding half-twisted bands, the ordering of the vertices is to be read from bottom to top, and the ordering of the edges from left to right. Right: an embedding of the spine as a strong deformation retract of the Seifert surface.}
    \label{fig.1colorspine}
\end{figure}

The resulting Seifert surface~$S$ can be encoded as a graph~$G$ with $n$ vertices, each corresponding to a Seifert disk, and one edge for each half-twisted band. We call this graph the \textbf{spine} of the Seifert surface. We can fully encode $S$ by decorating each edge of the spine with either a plus sign or a minus sign to record the handedness of the corresponding half-twist, and by remembering the vertical ordering of the Seifert disks and the ordering of the edges around the stack of disks. See Figure~\ref{fig.1colorspine} (left) for an example. We will refer to this package of data as the \textbf{decorated spine} for the surface. We also remark that the spine embeds naturally as a strong deformation retract of the surface, as illustrated in Figure~\ref{fig.1colorspine} (right). We draw the embedding in such a way that the vertices of the spine all lie to the right of the edges. Later, it will turn out to be convenient that we have adopted one such choice.

\subsection{Constructing C-complexes}\label{sec.fullalg}

\begin{figure}[h]
 \centering
\def\svgwidth{0.5\linewidth}
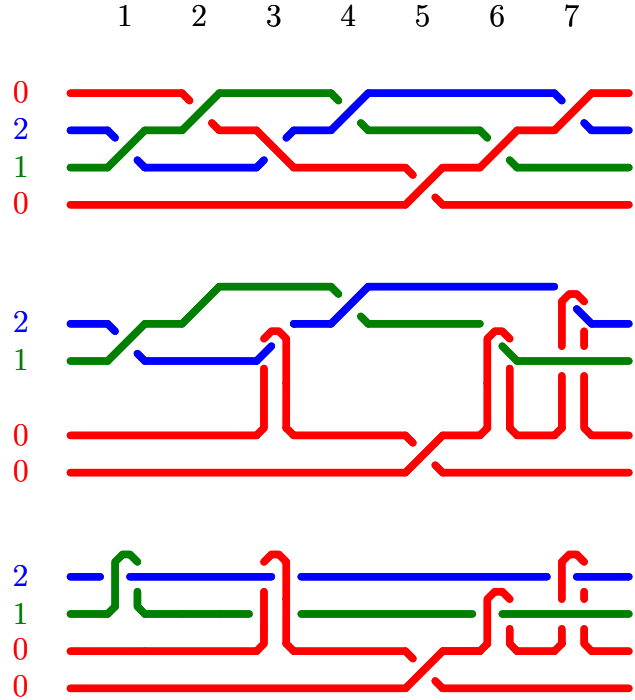
 \caption{Sorting the strands by color. Top: a braid on $n=4$ strands and $\mu = 3$ colors, with crossing sequence $(-2, -3, 2, -3, -1, -2, -3)$. Numbers on the left indicate the color of the strand starting on that position. Middle: dragging the $0$-colored strands across the back produces hooks at the points where they crossed over strands of different color. Bottom: repeating the procedure for all colors yields a diagram with crossings between strands of the same color and hooks between strands of different colors. Here, crossing~$\#5$ is preserved, crossings $\#2$ and $\#4$ do not show up in the final diagram because the higher-colored strand crosses over a lower-colored strand, crossing~$\#3$ produces a right-handed hook, and crossings $\#1$,~$\#6$~and~$\#7$ produce left-handed hooks. In the original braid, crossing $\#3$ occurs below the $1$-colored strand, so the resulting hook appears in front of that strand. In contrast, crossing~$\#7$ occurs above the $1$-colored strand, and hence gives rise to a hook passing behind that strand.}
\label{fig.multicolorbraid}
\end{figure}

We now explain how to generalize the previous con\-stru\-ction to colored links. This time, besides the input data of the number~$n$ of strands in the braid and the sequence $s = (s_1, \ldots, s_m)$ of crossings, we also have the data of a $\mu$-coloring of the braid. What this means in practice is that if $\sigma \in \Sigma_n$ is the permutation induced by the braid, then to each orbit of $\sigma$ we associate a color in~$\{0, \ldots, \mu-1\}$. See Figure~\ref{fig.multicolorbraid} (top) for an example.

\subsubsection{From crossings to hooks}\label{sec.dragdown}

To understand the translation of the input data into a decorated graph, it is useful to draw the braid diagram with the strands sorted by color, as we now explain. We start by considering all the strands with color~$0$, and isotope them as a stack to the bottom of the braid, keeping everything else fixed in place. We adopt the convention that the strands are to be moved \emph{across the back} of the braid. This is not an isotopy relative endpoints of the braid, but the endpoints are moved in parallel, so it extends to an isotopy of the braid closure. This modification does not affect crossings between $0$-colored strands, however some of these strands might get caught in differently-colored strands, leaving hook formations where some of the crossings used to be -- see Figure \ref{fig.multicolorbraid} (middle). Specifically, each time a $0$-colored strand crosses \emph{over} a strand with a different color, that crossing will appear as a hook in the modified diagram. The handedness of this hook depends on whether the $0$-colored strand was moving one position up or down. On the other hand, crossings of $0$-colored strands \emph{under} strands of different colors will not manifest in the final picture.

Having moved the $0$-colored strands to the bottom of the picture, we then proceed by moving the $1$-colored strands, as a stack, to the space above the $0$-colored strands, and so on until all colors have been moved. In the end, we obtain a diagram for a braid whose closure is isotopic to that of our starting braid, and where each position contains only strands of a fixed color -- see Figure~\ref{fig.multicolorbraid} (bottom). Moreover, the only interactions between different strands are either crossings between strands of the same color (which may be positive or negative), or hooks of some strand around a strand of a higher color. These hooks have a handedness, as already explained, may span several positions in the braid, and may travel across the front or the back of the strands in intermediate positions. This last distinction is a reflection of whether the crossing that originated the hook occurred above or below the intermediate strand.

\subsubsection{Filling in the surfaces}\label{sec.fillin}

\begin{figure}[h]
    \centering
    \def\svgwidth{0.60\linewidth}
    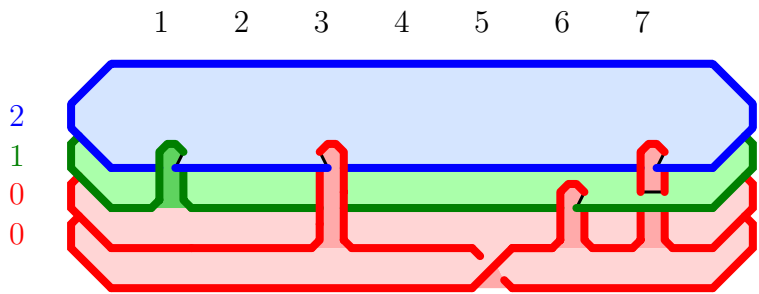
    \caption{Constructing colored surfaces from a braid with strands sorted by color. Following up the example from Figure~\ref{fig.multicolorbraid}, we see a Seifert disk for each position in the braid, each now having a well-defined color. As in the single-color case, crossing~$\#5$ gives rise to a half-twisted band between two disks of the same color. The hooks from crossings $\#1, \#3, \#6$~and~$\#7$ yield fingers connecting disks of different colors, with a clasp intersection at the apex. Moreover, since the hook from crossing $\#7$ passes behind the $1$-colored strand, it also produces a ribbon intersection between the finger and the $1$-colored disk. The fingers arising from crossings~$\#1, \#3$~and~$\#6$ are clean.}
    \label{fig.colorpancakes}
\end{figure}

In a similar fashion to what was done in the single-color case, we can now stare at such a diagram and visualize a collection of (possibly disconnected) orientable surfaces bounded by link components of the same color, with surfaces of different colors intersecting only along clasps or ribbons, as in Figure~\ref{fig.colorpancakes}. More explicitly, these surfaces are constructed by starting with a disk for each of the $n$ positions, filling in the crossings between strands of the same color with half-twisted bands, and filling in the hooks with protrusions of the disks, which we will call \textbf{fingers}. Each finger forms a clasp with the Seifert disk at the position where its boundary hooks, and whenever the finger passes behind some strand, it creates a ribbon intersection with the corresponding disk. A finger without ribbon intersections is said to be \textbf{clean}.

\subsubsection{Removing ribbon intersections}

In order to turn this collection of surfaces into a C-complex, we need to exchange the ribbon intersections that show up by clasp intersections. It turns out that the framework developed so far can neatly handle this situation, because of the following observation:

\begin{obs*}
Let $1 \le k_1 < k_2 < k_3 \le n$. Suppose a finger from disk $\# k_1$ whose boundary hooks with disk~$\# k_3$ has its bottom-most ribbon intersection with disk $\# k_2$. Then the isotopy type of the link does not change when that ribbon intersection is removed, and two clean fingers from disk $\# k_1$ to disk~$\# k_2$ are added: a right-handed one to the left of the original finger, and a left-handed one to the right. 
\end{obs*}

\begin{figure}[h]
    \centering
        \def\svgwidth{0.6\linewidth}
    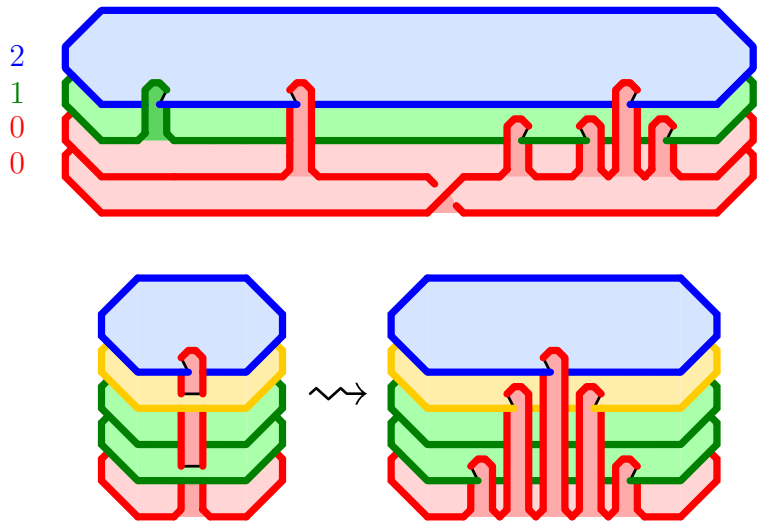
    \caption{Exchanging ribbon intersections for clean fingers. Top: in the running example from Figure~\ref{fig.colorpancakes}, the only ribbon intersection occurs at the right-most finger. This ribbon intersection is removed at the cost of adding two clasp intersections given by clean fingers. Bottom: a more complicated example, where this trick is iterated over multiple ribbon intersections in a finger, from bottom to top.}
    \label{fig.cleanfingers}
\end{figure}

As illustrated in Figure~\ref{fig.cleanfingers}, by sequentially applying this observation to the ribbon inter\-sections in our colored surfaces from the bottom to the top, we reach a setting where all fingers are clean, and so there are only clasp intersections between the surfaces. 

\subsubsection{Cleaning up}

\begin{figure}[h]
    \centering
        \def\svgwidth{0.4\linewidth}
\begingroup%
  \makeatletter%
  \providecommand\color[2][]{%
    \errmessage{(Inkscape) Color is used for the text in Inkscape, but the package 'color.sty' is not loaded}%
    \renewcommand\color[2][]{}%
  }%
  \providecommand\transparent[1]{%
    \errmessage{(Inkscape) Transparency is used (non-zero) for the text in Inkscape, but the package 'transparent.sty' is not loaded}%
    \renewcommand\transparent[1]{}%
  }%
  \providecommand\rotatebox[2]{#2}%
  \newcommand*\fsize{\dimexpr\f@size pt\relax}%
  \newcommand*\lineheight[1]{\fontsize{\fsize}{#1\fsize}\selectfont}%
  \ifx\svgwidth\undefined%
    \setlength{\unitlength}{54.70495569bp}%
    \ifx\svgscale\undefined%
      \relax%
    \else%
      \setlength{\unitlength}{\unitlength * \real{\svgscale}}%
    \fi%
  \else%
    \setlength{\unitlength}{\svgwidth}%
  \fi%
  \global\let\svgwidth\undefined%
  \global\let\svgscale\undefined%
  \makeatother%
  \begin{picture}(1,0.39758962)%
    \lineheight{1}%
    \setlength\tabcolsep{0pt}%
    \put(0,0){\includegraphics[width=\unitlength,page=1]{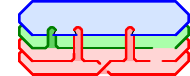}}%
    \put(-0.00385826,0.08241414){\color[rgb]{1,0,0}\makebox(0,0)[lt]{\lineheight{1.25}\smash{\begin{tabular}[t]{l}$0$\end{tabular}}}}%
    \put(-0.00385826,0.21525288){\color[rgb]{0,0.50196078,0}\makebox(0,0)[lt]{\lineheight{1.25}\smash{\begin{tabular}[t]{l}$1$\end{tabular}}}}%
    \put(-0.00385826,0.28463684){\color[rgb]{0,0,1}\makebox(0,0)[lt]{\lineheight{1.25}\smash{\begin{tabular}[t]{l}$2$\end{tabular}}}}%
    \put(-0.00385826,0.1501306){\color[rgb]{1,0,0}\makebox(0,0)[lt]{\lineheight{1.25}\smash{\begin{tabular}[t]{l}$0$\end{tabular}}}}%
    \put(0,0){\includegraphics[width=\unitlength,page=2]{clasp_cleanup.pdf}}%
  \end{picture}%
\endgroup%

    \caption{The result of removing the redundant pair of fingers from the surfaces in Figure~\ref{fig.cleanfingers} (top). In this example, this completes the cleanup step.}
    \label{fig.clasp_cleanup}
\end{figure}

At this point, we might not yet be in the presence of a C-complex, because the surfaces in each color might fail to be connected, and thus fail to be Seifert surfaces. Moreover, the union of all surfaces might not be connected. 
Before addressing that, however, we perform a clean-up step that is relevant only for reasons of computational efficiency: we simplify our surfaces by removing consecutive pairs of oppositely oriented half-twisted bands or fingers between the same disks (say, by scanning from left to right).
This step is carried out cyclically, that is, if the first and last bands or fingers connect the same pair of disks and have opposite signs, they are also removed. This step is repeated until no such redundant pairs remain. Clearly this does not change the isotopy type of the link. We exemplify on our running example in Figure~\ref{fig.clasp_cleanup}.

\subsubsection{Guaranteeing connectedness conditions}

Next, we make sure that all Seifert disks of the same color are connected by at least one half-twisted band. This is done exactly as in the single color case: whenever two adjacent disks of the same color are not connected to one another, add a pair of half-twisted bands with opposite handedness between them. In this way, we end up with a Seifert surface for each color.

\begin{figure}[h]
    \centering
    \def\svgwidth{0.8\linewidth}
    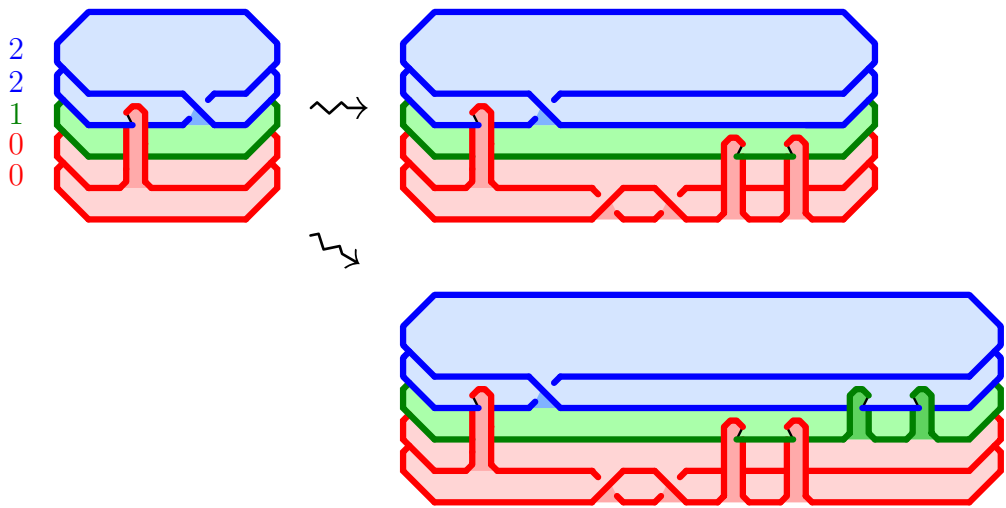
    \caption{Ensuring connectedness conditions. Top: the $0$-colored surface on the left is not connected, so we turn it into a Seifert surface by adding a pair of oppositely-oriented half-twisted bands between the two $0$-colored disks. The union of the resulting Seifert surfaces then has two connected components, which we connect by adding a pair of consecutive oppositely oriented fingers between the disks in positions $1$~and~$3$.
    Bottom: if we insist that all Seifert surfaces in the C-complex have non-empty intersection, we should also introduce a pair of fingers between the disks in positions $3$~and~$4$.}
    \label{fig.connectedcolors}
\end{figure}

The next step is to ensure that the C-complex itself is connected. To this end, we sort the Seifert surfaces into the connected components of the C-complex, and add a pair of consecutive oppositely oriented clean fingers between the bottom-most disk of the bottom-most component, and the bottom-most disk of each of the other components. As before, this does not change the isotopy type of the link. The result of the processing form this and the previous paragraph is exemplified with a toy example in Figure~\ref{fig.connectedcolors} (top).

As mentioned at the end of Section~\ref{sec.invariants}, one might wish to use the generalized Seifert matrices for computing a presentation matrix for the multivariable Alexander module $\h_1(S^3\setminus L;\Lambda_\mu)$ \cite[Theorem~3.2]{CF08}. This application, however, requires not only that the C-complex be connected, but in fact that every two Seifert surfaces have non-empty intersection. If this stronger condition is desired, our algorithm introduces, additionally, pairs of fingers between the bottom-most disks of every two disjoint Seifert surfaces, as exemplified in Figure~\ref{fig.connectedcolors} (bottom).

\subsubsection{Encoding the C-complex as a decorated spine}

\begin{figure}[h]
    \centering
    \def\svgwidth{0.6\linewidth}
    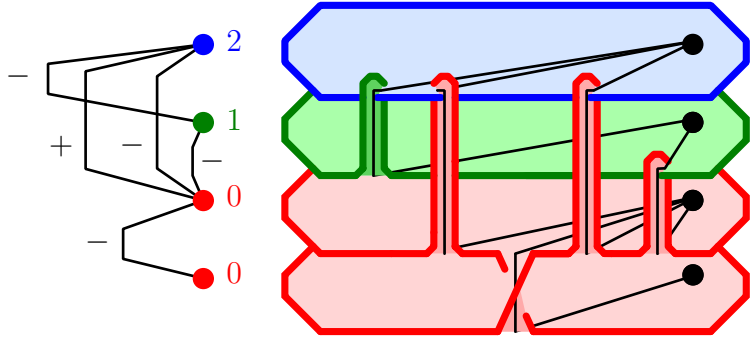
    \caption{A decorated spine for a C-complex. Left: the decorated spine for the C-complex in Figure~\ref{fig.clasp_cleanup}.  Vertices are to be read as ordered from bottom to top, and edges from left to right. This time, vertices are labeled with colors, which we indicate to the right. Edges represent half-twisted bands or fingers depending on whether their endpoints are the same color. Signs indicate the handedness of the half-twisted band or finger. Right: an embedding of the spine in the C-complex as a strong deformation retract.}
    \label{fig.multicolorspine}
\end{figure}

We are now ready to encode our C-complex as a decorated graph, which we do by extending the definition of a decorated spine to the multi-colored setting. We have one vertex for each Seifert disk, remembering their bottom-to-top ordering, as well as their color. Now we add one edge for each half-twisted band or finger, remembering the left-to-right order. We don't have to explicitly remember whether each edge represents a half-twisted band or a finger, as that is determined by whether its endpoints are of the same color. We do however need to store a sign for each edge, as before, to encode the handedness of the corresponding half-twisted band or finger. In Figure~\ref{fig.multicolorspine} we give an example of a decorated spine, and also illustrate the fact that the spine embeds as a strong deformation retract of the clasp complex. In the schematic, we again drew the vertices of the embedded spine to the right of all edges, for reasons that will become apparent in the next subsection. 

\subsection{Reading off a Seifert matrix}\label{sec.readmatrix}

We now explain how to read off a Seifert matrix for the C-complex~$S$ from its decorated spine. 

First, one must choose a basis for the first homology~$\h_1(S)$. Since $S$~strongly deformation  retracts onto its spine~$G$, this amounts to finding a basis for the homology of a finite connected graph. This is a routine exercise, but we give a quick sketch: choose a maximal tree~$T$ for~$G$, and then $G/T$~is a wedge of, say,  $r$~circles, with the collapse map $G \to G/T$ being a homotopy equivalence. After orienting each of the $r$~circles, the resulting collection of oriented loops gives a basis of~$\h_1(G/T)$. Now, each loop lifts to an edge~$e$ of~$G$ that is not in~$T$, say oriented from the vertex~$v_1$ to the vertex~$v_2$. The corresponding element of~$\h_1(G)$ is represented by the circuit obtained by concatenating~$e$ with the unique path in~$T$ from~$v_2$ to~$v_1$.

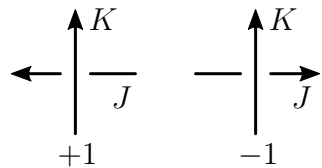
\begin{figure}[h]
    \centering
    \captionsetup{margin=20pt}
    \def \svgwidth{0.25 \linewidth}
\begingroup%
  \makeatletter%
  \providecommand\color[2][]{%
    \errmessage{(Inkscape) Color is used for the text in Inkscape, but the package 'color.sty' is not loaded}%
    \renewcommand\color[2][]{}%
  }%
  \providecommand\transparent[1]{%
    \errmessage{(Inkscape) Transparency is used (non-zero) for the text in Inkscape, but the package 'transparent.sty' is not loaded}%
    \renewcommand\transparent[1]{}%
  }%
  \providecommand\rotatebox[2]{#2}%
  \newcommand*\fsize{\dimexpr\f@size pt\relax}%
  \newcommand*\lineheight[1]{\fontsize{\fsize}{#1\fsize}\selectfont}%
  \ifx\svgwidth\undefined%
    \setlength{\unitlength}{77.51919604bp}%
    \ifx\svgscale\undefined%
      \relax%
    \else%
      \setlength{\unitlength}{\unitlength * \real{\svgscale}}%
    \fi%
  \else%
    \setlength{\unitlength}{\svgwidth}%
  \fi%
  \global\let\svgwidth\undefined%
  \global\let\svgscale\undefined%
  \makeatother%
  \begin{picture}(1,0.50879449)%
    \lineheight{1}%
    \setlength\tabcolsep{0pt}%
    \put(0,0){\includegraphics[width=\unitlength,page=1]{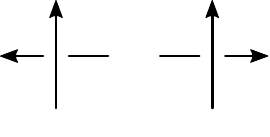}}%
    \put(0.83172968,0.44304284){\color[rgb]{0,0,0}\makebox(0,0)[lt]{\lineheight{1.25}\smash{\begin{tabular}[t]{l}$K$\end{tabular}}}}%
    \put(0.90461263,0.18939041){\color[rgb]{0,0,0}\makebox(0,0)[lt]{\lineheight{1.25}\smash{\begin{tabular}[t]{l}$J$\end{tabular}}}}%
    \put(0.25122942,0.44304284){\color[rgb]{0,0,0}\makebox(0,0)[lt]{\lineheight{1.25}\smash{\begin{tabular}[t]{l}$K$\end{tabular}}}}%
    \put(0.32411627,0.18939041){\color[rgb]{0,0,0}\makebox(0,0)[lt]{\lineheight{1.25}\smash{\begin{tabular}[t]{l}$J$\end{tabular}}}}%
    \put(0.73190411,0.00944449){\color[rgb]{0,0,0}\makebox(0,0)[lt]{\lineheight{1.25}\smash{\begin{tabular}[t]{l}$-1$\end{tabular}}}}%
    \put(0.15140218,0.00944449){\color[rgb]{0,0,0}\makebox(0,0)[lt]{\lineheight{1.25}\smash{\begin{tabular}[t]{l}$+1$\end{tabular}}}}%
  \end{picture}%
\endgroup%

    \caption{The contributions of the crossings in a diagram for the knots~$J, K$ to the linking number~$\lk(J, K)$.}
    \label{fig.crossingsign}
\end{figure}

We are now armed with a collection of circuits in the spine~$G$, and wish to read off the linking numbers of the corresponding embedded curves in the C-complex~$S$, and their push-offs. We remind the reader that the linking number~$\lk(J,K)$ of two knots~$J, K$ can be computed from a diagram as the signed count of the number of crossings of~$K$ over~$J$, with the sign convention depicted in Figure \ref{fig.crossingsign}. We will use diagrams for the curves in our homology basis for~$S$ and their push-offs obtained from depictions as in Figure~\ref{fig.multicolorspine}. The local nature of the linking number computation is then well-suited to our description of~$\h_1(S)$ as a collection of circuits in~$G$, because we can simply count the contribution of each pair of edges. More explicitly, let $\gamma, \delta$~be circuits in~$G$ given as sequences of oriented edges $\gamma = (\alpha_{1}, \ldots, \alpha_{r}), \delta = (\beta_1, \ldots, \beta_{s})$. For a $\mu$-tuple of signs $\epsilon \in \{\pm 1\}^\mu$, we have
\[ \lk(\gamma^\epsilon, \delta) = \sum_{\substack{1 \le i \le r \\ 1 \le j \le s}} [ \alpha_i^\epsilon, \beta_j ], \label{linknumberformula} \] 
where the $\epsilon$-superscript denotes the push-off dictated by~$\epsilon$, as described in Section~\ref{sec.seifertmatrices}, and for oriented edges~$\alpha, \beta$ of~$G$, the symbol~$[\alpha^\epsilon, \beta]$ denotes the signed number of crossings of~$\beta$ over~$\alpha^\epsilon$ in some fixed diagram for $\gamma$~and~$\delta$.

At this point, we have reduced the computation of the desired linking numbers $\lk(\gamma^\epsilon, \delta)$ to counting the signed crossings~$[\alpha^\epsilon, \beta]$ for oriented edges~$\alpha, \beta$ of~$G$. We emphasise that the symbol~$[\alpha^\epsilon, \beta]$ makes sense only under a choice of diagram, which should be globally fixed in order for the above formula for linking numbers to hold. This is where it is relevant that we insist on always drawing the vertices of the spine~$G$ to the right of the edges, as in Figures~\ref{fig.1colorspine}~and~\ref{fig.multicolorspine}. Having fixed this convention, we proceed to investigate how to read off the symbols~$[\alpha^\epsilon, \beta]$ from $G$~and its decorations.

Denote the endpoints of~$\alpha$ by $u_1, u_2$, with $\alpha$ oriented from~$u_1$ to~$u_2$, and similarly suppose $\beta$~goes from~$v_1$ to~$v_2$. We assume that $\alpha$~and~$\beta$ are both oriented ``upwards'' in~$G$, that is, for the total vertex order packaged into the decoration of~$G$, we have $u_1 < u_2$~and~$v_1 < v_2$. The other cases are recovered from the obvious identities $[\overline{\alpha}^\epsilon, \beta] = [\alpha^\epsilon, \overline{\beta}] = - [\alpha^\epsilon, \beta]$, where over-lines indicate orientation reversal.

Observe now that in our picture of the spine~$G$ embedded in~$S$, the four points in the set $P:=\{u_1^\epsilon, u_2^\epsilon, v_1, v_2\}$ are all distinct and vertically aligned. This is true even if one of the~$u_i$ equals one of the~$v_j$, because then the push-off given by the sign~$\epsilon_i$ moves~$u_i$ above or below~$v_j$. Reading these points from bottom to top, one thus obtains a total order on~$P$ that is easily recovered from the decoration of~$G$: we first order~$P$ partially by reading the indices of $u_1, u_2, v_1, v_2$, and then break potential ties, necessarily between a~$u_i$ and a~$v_j$, by reading the sign~$\epsilon_i$.

One now sees that $[\alpha^\epsilon, \beta]$~can only be non-zero if the~$u_i^\epsilon$ and the~$v_j$ are alternating, that is, if $u_1^\epsilon < v_1 < u_2^\epsilon < v_2 $ or $v_1 < u_1^\epsilon < v_2 < u_2^\epsilon$. Let us thus assume we are in this situation, and consider first the case where the edges $\alpha$~and~$\beta$ are distinct. In this case, it is clear that one among~$\alpha^\epsilon, \beta$ crosses over the other precisely once. The value of~$[\alpha^\epsilon, \beta]$ is then non-zero exactly if it is~$\beta$ crossing over~$\alpha^\epsilon$. Now, our convention of drawing the vertices of~$G$ to the right of all edges implies that $\beta$~crosses over~$\alpha^\epsilon$ precisely if in the total ordering of edges of~$G$ we have~$\alpha < \beta$. We then see by direct inspection that 

\[ [\alpha^\epsilon, \beta] = \begin{cases} 1 &\text{if $ v_1 < u_1^\epsilon < v_2 < u_2^\epsilon$,} \\ -1 &\text{if $ u_1^\epsilon< v_1 < u_2^\epsilon < v_2$.} 
\end{cases}\]

\begin{figure}[h]
    \centering
    \captionsetup{margin=20pt}
    \def \svgwidth{0.55 \linewidth}
    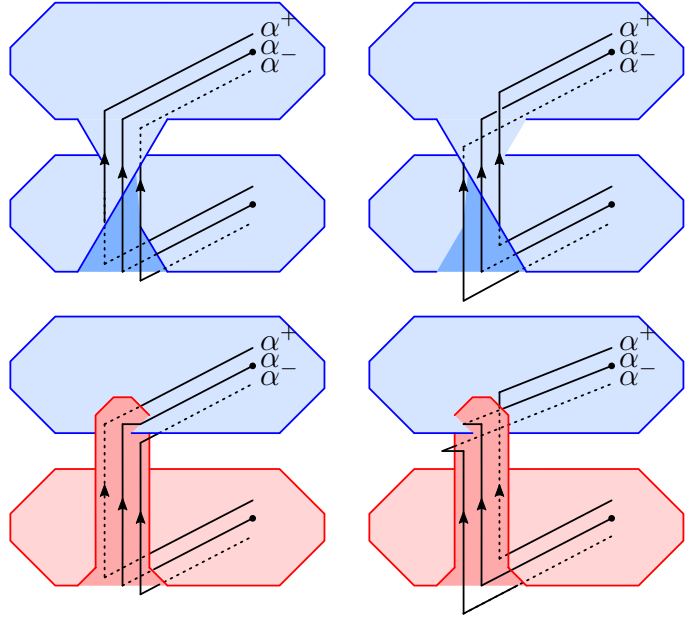
    \caption{An upwards-oriented edge~$\alpha$ of the spine~$G$, together with its push-offs $\alpha^+, \alpha^-$. We consider the cases where $\alpha$~corresponds to a half-twisted band (top) or a finger (bottom), and whether this band/finger is left-handed (left) or right-handed (right).}
    \label{fig.pushoffs}
\end{figure}

If~$\alpha = \beta$, we must analyse several cases. Still under the assumption that $\alpha$~is oriented upwards and that $\alpha, \alpha^\epsilon$~have endpoints that alternate along the vertical direction, we have to consider:
\begin{itemize}
    \item whether $\alpha$~corresponds to a half-twisted band or a finger,
    \item the handedness of the half-twisted hand or finger,
    \item whether the sign of $\epsilon$~at the endpoints of~$\alpha$ is $+$~or~$-$ (in the case of fingers, the assumption that the endpoints of $\alpha,\alpha^\epsilon$ alternate implies that the push-off direction dictated by~$\epsilon$ is the same at both endpoints).
\end{itemize}
This amounts to eight cases, which are depicted in Figure~\ref{fig.pushoffs}. By direct inspection, we obtain the results in the following table:

\begin{center}
\begin{tabular}{r|l | l}
& left-handed & right-handed\\ \hline
\multirow{2}{*}{half-twisted band} & $[\alpha^+, \alpha] = +1$ & $[\alpha^+, \alpha] = 0$\\
& $[\alpha^-, \alpha] = 0$ & $[\alpha^-, \alpha] = -1$\\ \hline
\multirow{2}{*}{finger} & $[\alpha^+, \alpha] = +1$ & $[\alpha^+, \alpha] = 0$\\
& $[\alpha^-, \alpha] = 0$ & $[\alpha^-, \alpha] = -1$\\\end{tabular}
\end{center}

With this, we finish the explanation of how to determine the value of~$[\alpha^\epsilon, \beta]$ (in a diagram following our conventions) for any two oriented edges in~$G$, just by reading the decoration of $G$. Hence, by the linking number formula on page \pageref{linknumberformula},  we know how to compute $\lk(\gamma^\epsilon, \delta)$ for any circuits $\gamma, \delta$ in~$G$. Applying this to a homology basis for~$G$ we produce the desired generalized Seifert matrix.

\section{Additional comments on the implementation}\label{section:implementation}

In this brief section we say a few words about the actual computer implementation of our algorithm. 

\subsection{Input format}
In \clasper, the input format for the braids follows the convention of the ``braid notation'' in LinkInfo \cite{LM22} and of the website ``Seifert Matrix Computations'' (SMC) \cite{CKL16}. Note that in the explanation of the notation in SMC, the positions in the braid are numbered from top to bottom, but the sign convention for left/right-handed crossings is the same as ours. Hence, given a sequence of crossings, the braid specified by our convention and the one specified as in SMC differ merely by a rotation of half a turn about a horizontal line in the projection plane, which is immaterial.

\subsection{Output format}
\clasper\ displays the colored link invariants on the graphical interface, and also allows the user to them as \LaTeX\ code.

The button ``Export Seifert matrices'' allows the user to save a text file containing a presentation matrix for the multivariable Alexander module $\h_1(S^3\setminus L;\Lambda_\mu)$ (where $\Lambda_\mu:=\Z[t_1^{\pm 1},\dots,t_\mu^{\pm 1},(1-t_1)^{-1},\dots,(1-t_\mu)^{-1}]$), and the collection of generalized Seifert matrices used to compute it, in a format compatible with SageMath.
Above each generalized Seifert matrix is indicated the sign tuple $\epsilon \in \{0,1\}^\mu$ to which it corresponds. For the running example in this section, the output looks as follows.

\bigskip
\begin{quote}
\begin{footnotesize}
\begin{verbatim}
Presentation Matrix
Matrix([[0, t0*t1*t2 - t0*t2], [t1 - 1, -t0*t1*t2 + 1]])


Generalized Seifert Matrices

[-1, -1, -1]
Matrix([[0, -1], [0, 1]])

[-1, -1, 1]
Matrix([[0, 0], [0, 0]])

[-1, 1, -1]
Matrix([[0, -1], [0, 0]])

[-1, 1, 1]
Matrix([[0, 0], [0, 0]])

[1, -1, -1]
Matrix([[0, 0], [0, 0]])

[1, -1, 1]
Matrix([[0, 0], [-1, 0]])

[1, 1, -1]
Matrix([[0, 0], [0, 0]])

[1, 1, 1]
Matrix([[0, 0], [-1, 1]])
\end{verbatim}
\end{footnotesize}
\end{quote}

\subsection{Optimizing determinant computations}
In our approach for determining the Conway potential function of a braid closure, the most computationally demanding step is the calculation of the determinant
\[\operatorname{det}\Big(-\sum_{\epsilon \in \{\pm 1\} ^\mu} A^\epsilon\cdot \epsilon_1 \dots \epsilon_\mu\cdot  t_1^{\epsilon_1}\dots t_\mu^{\epsilon_\mu}\Big),\]
which appears in the formula in Theorem~\ref{thm:conway-potential-function}. We employ the Bareiss algorithm for efficient computation of determinants using integer arithmetic.

Moreover, in order to try and streamline this step, \clasper\ computes several different spines for the braid closure, obtained by randomly permuting the colors of the link. In other words, when performing the step described in Subsection~\ref{sec.dragdown}, we ``drag down'' the colors in different orders. The idea is to find C-complexes whose spines have homology with small rank, so that the determinant computation is performed on smaller matrices. We do this by trying out 500 randomly chosen permutations of the colors, and selecting one with minimal homology rank.

\subsection{Signature computations and floating point arithmetic}
The computation of Cimasoni-Florens signatures is carried out in floating point arithmetic. \clasper\ will consider to be~$0$ any eigenvalue of absolute value below $10^{-5}$, but will also display the computed eigenvalues of the relevant matrix~$H(\omega)$ from Definition~\ref{dfn.CFsignatures}.

\subsection{Libraries used and download location}
\clasper\ was written by the second author in Python 3 using the libraries numpy, matplotlib, tkinter and sympy.
 A Windows installer and the Python source code are available at \url{https://github.com/Chinmaya-Kausik/py_knots}.


\begin{thebibliography}{10}
\bibitem[APS75]{APS75}
M. Atiyah, V. Patodi and I. Singer. {\em
Spectral asymmetry and Riemannian geometry. II}, 
Math. Proc. Camb. Philos. Soc. 78 (1975), 405--432.
\bibitem[CG75]{CG75} 
A. Casson and C. M. Gordon. {\em Cobordism of classical knots in $S^3$}, Printed notes. Orsay (1975).
\bibitem[CG78]{CG78}
A. Casson and C. M. Gordon. {\em On slice knots in dimension three}, Proc. Symp. in Pure Math.
XXX (1978), Part 2, 39–53. 
\bibitem[BZH14]{BZH14}
G. Burde, H. Zieschang and M. Heusener. {\em 
Knots}, 3rd fully revised and extended edition.
De Gruyter Studies in Mathematics 5. Berlin: Walter de Gruyter (2014). 
\bibitem[Ci04]{Ci04}
D. Cimasoni, {\em A geometric construction of the Conway potential function}, Comment. Math.
Helv. 79 (2004), 
 124–146.
\bibitem[CF08]{CF08}
D. Cimasoni and V. Florens.
{\em Generalized Seifert surfaces and signatures of colored links},
Trans. Am. Math. Soc. 360, No. 3, 1223-1264 (2008). 
\bibitem[Col16]{Col16}
J. Collins. {\em 
An algorithm for computing the Seifert matrix of a link from a braid representation},
Ghys (ed.) et al., Six papers on signatures, braids and Seifert surfaces. Rio de Janeiro: Sociedade Brasileira de Matem\'atica (SBM), Ensaios Matem\'aticos 30 (16), 246-262.
\bibitem[CKL16]{CKL16}
J. Collins, T. K\"oppe and L. Lewark.
{\em Seifert Matrix Computations}.\\
\url{https://www.maths.ed.ac.uk/~v1ranick/julia/index.htm}
\bibitem[Con18]{Con18}
A. Conway. {\em 
An explicit computation of the Blanchfield pairing for arbitrary links},  Can. J. Math. 70 (2018), 983-1007. 
\bibitem[CFT18]{CFT18}
A. Conway, S. Friedl and E. Toffoli. {\em 
The Blanchfield pairing of colored links},  Indiana Univ. Math. J. 67 (2018),  2151-2180. 
\bibitem[Con70]{Con70}
J. H. Conway. {\em An enumeration of knots and links, and some of their algebraic properties},
in: Computational Problems in Abstract Algebra (Proc. Conf., Oxford, 1967), 329--358,
Pergamon, Oxford, 1970.
\bibitem[Coo82]{Coo82}
D. Cooper. {\em  The universal abelian cover of a link}, Low-dimensional topology (Bangor, 1979), London Math. Soc. Lecture Note Ser., 48, Cambridge Univ. Press (1982), 51-66.
\bibitem[DMO21]{DMO21}
C. Davis, T. Martin and C. Otto. {\em Moves relating C-complexes: A correction to Cimasoni's ``A geometric construction of the Conway potential function"}, Preprint (2021),
arXiv:2105.10495.
\bibitem[Ka96]{Ka96}
A. Kawauchi. {\em A survey of knot theory}, Birkh\"auser (1996).
\bibitem[Li97]{Li97}
W. B. R. Lickorish. {\em  An introduction to knot theory}, Springer Graduate Text in Mathematics (1997). 
\bibitem[LM22]{LM22}
C. Livingston and A. H. Moore. {\em LinkInfo: Table of Link Invariants}, \url{https://linkinfo.sitehost.iu.edu} (April 14, 2022). 
\bibitem[O'B02]{O'B02}
K. O'Brien. {\em Seifert's Algorithm, Ch\^atelet Bases and the Alexander Ideals of Classical Knots}, PhD thesis, University of Durham (2002).\\
\url{http://etheses.dur.ac.uk/4192/}
\bibitem[Ro90]{Ro90}
D. Rolfsen. {\em 
Knots and links},  2nd print. with corr.
Mathematics Lecture Series. 7. Houston, TX: Publish or Perish. xiv (1990). 
\bibitem[Se34]{Se34}
H. Seifert. {\em \"Uber das Geschlecht von Knoten},  Math. Ann. 110 (1934), 571-592.

\end{thebibliography}
\end{document}